\documentclass[12pt]{amsart}

\usepackage {hyperref}

\usepackage{amssymb}

\usepackage{enumerate}

\makeatletter
\@namedef{subjclassname@2010}{%
  \textup{2010} Mathematics Subject Classification}
\makeatother

\newtheorem {theorem} {Theorem}

\newtheorem {proposition} [theorem] {Proposition}
\newtheorem {definition} [theorem] {Definition}

\newcommand {\apclass} [1] {\ensuremath{\mathrm A_{#1}}}
\newcommand {\frclass} [2] {\ensuremath{\mathrm F^{#1}_{#2}}}

\newcommand {\lclass} [2] {\ensuremath{\mathrm L_{#1} \left( #2 \right) }}

\newcommand {\lclassg} [1] {\ensuremath{\mathrm L_{#1}}}

\newcommand {\BMO} {\ensuremath {\mathrm {BMO}}}

\DeclareMathOperator* {\esssup} {ess\,sup}
\DeclareMathOperator* {\supp} {supp}

\newcommand {\weightu} {\ensuremath {\mathit u}}
\newcommand {\weightv} {\ensuremath {\mathit v}}
\newcommand {\weightw} {\ensuremath {\mathit w}}

\newcommand {\abp} [2] {\ensuremath {B\apclass {#1}^{(\mathrm {MC})} \left(#2\right)}}
\newcommand {\abpo} [2] {\ensuremath {B\apclass {#1} \left(#2\right)}}
\newcommand {\fbab} [3] {\ensuremath {B\frclass {#1} {#2} \left(#3\right)}}

\begin{document}


\baselineskip=17pt


\title [$\apclass {1}$-regularity and Calder\'on-Zygmund operators]
{$\mathbf A_1$-regularity and boundedness\\of Calder\'on-Zygmund operators. II}

\author [D.~V.~Rutsky] {Dmitry V. Rutsky}
\address {
Steklov Mathematical Institute\\
St. Petersburg Branch\\
Fontanka 27\\
191023 St. Petersburg, Russia}
\email {rutsky@pdmi.ras.ru}

\date{}

\begin{abstract}
Proof is given for the ``only if'' part of the result stated in the previous paper of the series
that a suitably nondegenerate Calde\-r\'on-Zygmund operator $T$ is bounded in a Banach lattice $X$
on $\mathbb R^n$ if and only if the Hardy-Littlewood maximal operator $M$ is bounded in both $X$ and $X'$,
under the assumption that $X$ has the Fatou property and $X$ is $p$-convex and $q$-concave with some $1 < p, q < \infty$.
We also get rid of an application of a fixed point theorem in the proof of the main lemma
and give an improved version of an earlier result concerning the divisibility of $\BMO$-regularity.
\end{abstract}

\subjclass[2010] {Primary 46B42, 42B25, 42B20, 46E30, 47B38}

\keywords {$\apclass {p}$-regularity, Calder\'on-Zygmund operator, Hardy-Littlewood maximal operator}

\maketitle

This paper is closely related to~\cite {rutsky2014sm}
and contains essentially no new non-technical results, hence
for the background and the generalities
we refer the reader to~\cite {rutsky2014sm}.

A.~Yu.~Karlovich and L.~Maligranda kindly pointed out to the author that
the proof of \cite [Theorem~16] {rutsky2014sm} has a flaw, namely that the relationship
$(X \lclassg {s})' = X' \lclassg {s'}$ is incorrect (and in fact it is always false).
Unfortunately, it is not clear if
\cite [Theorem~16] {rutsky2014sm} is true in the stated form.

Nevertheless, we will see that the main result of \cite {rutsky2014sm} is still true with only a slight
loss of generality concerning the nondegeneracy assumption imposed on a Calder\'on-Zygmund operator $T$.
Specifically, in place of $\apclass {2}$-nondegeneracy of $T$ (which is a condition that
the boundedness of $T$ in $\lclass {2} {\weightw^{-\frac 1 2}}$ implies $\weightw \in \apclass {2}$ with
an estimate for the constant)
we require that the kernels of both $T$ and its conjugate $T^*$
satisfy a standard assumption on growth along a certain singular direction
(see \cite [Chapter~5, \S 4.6] {stein1993}).
\begin {definition}
\label {nondegop}
We say that a singular integral operator $T$ on $\mathbb R^n$ is nondegenerate if
there exists a constant $c > 0$ and some $x_0 \in \mathbb R^n \setminus \{0\}$
such that for any ball $B \subset \mathbb R^n$ of radius $r > 0$ and any locally summable
nonnegative function $f$ supported on $B$
we have
\begin {equation}
\label {nondegsoeq}
|T f (x)| \geqslant c \frac 1 {|B|} \int_B f
\end {equation}
for all $x \in B \pm r x_0$.
\end {definition}
For example, the Hilbert transform $H$ and any of the Riesz transforms $R_j$ are 
nondegenerate in this sense.
A nondegenerate operator $T$ is also $\apclass {2}$-nondegenerate.  For details see \cite [Chapter~5, \S 4.6] {stein1993}.

\begin {theorem}
\label {themcr}
Suppose that $X$ is a Banach lattice of measurable functions on $\mathbb R^n \times \Omega$ that satisfies the Fatou property
and $X$ is $p$-convex and $q$-concave with some $1 < p, q < \infty$.
Let $T$ be a Calder\'on-Zygmund operator in $\lclass {2} {\mathbb R^n}$ such that
both $T$ and $T^*$ are nondegenerate.
The following conditions are equivalent.
\begin {enumerate}
\item
The Hardy-Littlewood maximal operator $M$ acts boundedly in $X$ and in the order dual $X'$ of $X$.
\item
All Calder\'on-Zygmund operators act boundedly in $X$.
\item
$T$ acts boundedly in $X$.
\end {enumerate}
\end {theorem}

Thus, concerning the necessity of $\apclass {1}$-regularity
we make no claims about the general spaces of homogeneous type, although in many cases a suitable generalization
of Definition~\ref {nondegsoeq} seems to be possible.
Another subtle loss of generality is that in contrast to~\cite [Theorem~16] {stein1993}
in the proof of $3 \Rightarrow 1$ we take advantage
of the assumption that $T$ is a Calder\'on-Zygmund operator as well as a nondegenerate operator, specifically that
$T$ is bounded in $\lclassg {t}$ for $1 < t < \infty$
with norm $O (t)$ as $t \to \infty$.

For the proof of $1 \Rightarrow 2 \Rightarrow 3$ of Theorem~\ref {themcr} see~\cite {rutsky2014sm}.
The proof of $3 \Rightarrow 1$ essentially follows the scheme of the flawed proof
of~\cite [Theorem~16] {rutsky2014sm}, but it seems to require a much more delicate approach
that we will present throughout the rest of the paper,
leading to the proof itself given at the end of
Section~\ref {noaro} below.
We briefly outline the structure of the argument, the details of which are also of some independent interest.

The following result was established (with some caveats) in~\cite [Theorem~A'] {rubiodefrancia1987};
a complete proof in the stated form can be found in \cite [Theorem~4] {rutsky2011en}.
Here and elsewhere $(S, \nu)$ is a space of homogeneous type and $(\Omega, \mu)$ is a $\sigma$-finite measurable space.
\begin {theorem}
\label {btsb}
Suppose that $Y$ is a Banach lattice on $(S \times \Omega, \nu \times \mu)$ with an order continuous norm.
If a linear operator $T$ is bounded in $Y^{\frac 1 2}$
then for every $f \in Y'$
there exists a majorant
$\weightw \geqslant |f|$, $\|\weightw\|_{Y'} \leqslant 2 \|f\|_{Y'}$, such that
$\|T\|_{\lclass {2} {\weightw^{-\frac 1 2}} \to \lclass {2} {\weightw^{-\frac 1 2}}} \leqslant C \|T\|_{Y^{\frac 1 2} \to Y^{\frac 1 2}}$,
where $C$ depends only on the Grothendieck constant $K_G$.
\end {theorem}

This yields almost at once the following version of Theorem~\ref {themcr} that we will need
in the proof of Theorem~\ref {themcr},
showing that Theorem~\ref {themcr}
is also valid for $p = 2$, $q = \infty$ (and by duality for $p = 1$ and $q = 2$), provided that
$X$ (respectively, $X'$) has order continuous norm.  The proof is given in Section~\ref {potmcr2} below.
\begin {theorem}
\label {themcr2}
Suppose that $X$ is a $2$-convex Banach lattice of measurable functions on  $(S \times \Omega, \nu \times \mu)$
having order continuous norm and the Fatou property.
Let $T$ be an $\apclass {2}$-non\-de\-gene\-rate linear operator in $\lclass {2} {S \times \Omega}$.
If $T$ acts boundedly in $X$ then the maximal operator $M$ is bounded in both~$X$ and~$X'$ with
a suitable estimate for the constants.
\end {theorem}

In contrast to \cite {rutsky2014sm}, in the present work we use the standard definition of 
the constant $[\weightw]_{\apclass {p}}$, $p > 1$ of a Muckenhoupt weight $\weightw \in \apclass {p}$
on $(S \times \Omega, \nu \times \mu)$ 
based on the Muckenhoupt condition:
$$
[\weightw]_{\apclass {p}} = \esssup_{\omega \in \Omega}
\sup_{B} \left(\frac 1 {\nu (B)} \int_B \weightw (\cdot, \omega) \right)
\left(\frac 1 {\nu (B)} \int_B \weightw^{-\frac 1 {p - 1}} (\cdot, \omega) \right)^{p - 1},
$$
where the supremum is taken over all balls $B \subset S$.

Recall that a quasi-normed lattice $X$ is called $\apclass {p}$ regular with constants $(C, m)$ if
every $f \in X$ admits a majorant $\weightw \in X$, $\weightw \geqslant |f|$, such that $\|\weightw\|_X \leqslant m \|f\|_X$
and
$\weightw$ belongs to the Muckenhoupt class $\apclass {p}$ with $[\weightw]_{\apclass {p}} \leqslant C$.

In Section~\ref {mlr}
we give (Proposition~\ref {a1apt}) a simplified proof of \cite [Proposition~8] {rutsky2014sm}
that does not use a fixed point theorem.
This yields a slightly improved version (Proposition~\ref {ainfainf}) of \cite [Proposition~12] {rutsky2014sm}
stating that
$\apclass {\infty}$-regularity of both $X$ and $X'$ implies
$\apclass {1}$-regularity of these lattices, where the assumption that $X$ satisfies the Fatou property is replaced
by a weaker assumption that $X'$ is a norming lattice for $X$.  Thus it suffices to establish
that condition~3 of Theorem~\ref {themcr} implies that $X'$ is $\apclass {\infty}$-regular;
interchanging $X$ with $X'$ would then show that $X$ is also $\apclass {\infty}$-regular.

Under condition~3 of Theorem~\ref {themcr2} we may apply Theorem~\ref {btsb} to lattice $Y = X^r \lclassg {s}$
with some fixed $r > 1$ sufficiently close to $1$ and all sufficiently large $s$, since $T$ is bounded in $Y^{\frac 1 2}$
by interpolation with some estimate for the norm that grows with $s$.
This yields $\apclass {2}$-regularity of $Y' = (X^{rs'})'^{\frac 1 {s'}}$, with an estimate on the growth of the constant
$C_s$
as $s \to \infty$.
Now the key idea is to show that the $\apclass {2}$-majorants $\weightw$ of functions from $Y'$ also satisfy the reverse
H\"older inequality with exponent $s'$ for some sufficiently large $s$,
which would yield $\apclass {2}$-regularity of $(X^u)'$, $u = r s'$,
and thus the required $\apclass {2}$-regularity of the lattice $X' = (X^u)'^{\frac 1 u} \lclassg {1}^{1 - \frac 1 u}$.

However, as discussed in Section~\ref {apono} below, in order 
to get an estimate for $C_s$ with a suitable rate of growth we also need to make sure that the weight $\weightw$
appearing in the conclusion of Theorem~\ref {btsb} (applied to $Y$) satisfies some additional assumptions, namely that
$\weightw^{-1}$
is a doubling weight with a constant independent of $s$.
Theorem~\ref {themcr2} allows us to obtain $\apclass {1}$-regularity
of $Y'^{\frac 1 2} \lclassg {t}^{\frac 1 2}$ from condition~3 of Theorem~\ref {themcr}
with a sufficiently large fixed $t$, where an estimate for the constants is independent of $s$.
An extension (Theorem~\ref {frdiv})
of \cite [Theorem~2] {rutsky2011en} concerning the divisibility of $\apclass {p}$-regularity,
which we introduce in Section~\ref {doapar} below,
allows us to prove that $Y'$ admits suitable majorants
$\weightw$ such that $\weightw^{-1} \in \apclass {3}$ 
(and hence $\weightw^{-1}$ is a doubling weight) with a constant independent of $s$,
and an adaptation (Theorem~\ref {btsbge}) of the original
fixed point argument from \cite [\S 2] {rutsky2014sm} makes it possible to impose this condition on the weights appearing in
the conclusion of Theorem~\ref {btsb}, thus completing the proof.

\section {Proof of Theorem~\ref {themcr2}}

\label {potmcr2}

Suppose that $X$ and $Y$ are normed lattices on a measurable space $\Omega$.
Lattice $Y$ is said to be \emph {norming} for $X$
if $f g \in \lclassg {1}$ for all $f \in X$ and $g \in Y$ and
$\|f\|_X = \sup_{g \in Y, \|g\|_Y = 1} \int_{\Omega} |f g|$ for all $f \in X$.
A normed lattice $X$ is always norming for its order dual $X'$.
Conversely, it is well known that $X'$ is a norming lattice for $X$
if $X$ satisfies either the Fatou property (implying that $X = X''$), or if $X$ is a Banach lattice
having order continuous norm (since then $X' = X^*$).
The fact that $\weightw \in \apclass {p}$ if and only if the maximal operator $M$ is bounded in
$\lclass {p} {\weightw^{-\frac 1 p}}$ with the appropriate estimates of the constants
yields at once the following result; see \cite [Proposition~13] {rutsky2011en}.
\begin {proposition}
\label {aptoconj}
Suppose that $X$ and $Y$ are normed lattices on $(S \times \Omega, \nu \times \mu)$
such that $Y$ is a norming space for $X$.
If $Y$ is $\apclass {p}$-regular with some $p > 1$ then $X^{\frac 1 p}$ is $\apclass {1}$-regular with
appropriate estimates for the constants.
\end {proposition}

The following result is a particular case of \cite [Proposition~13] {rutsky2015cen}; we give a complete proof for clarity.
\begin {proposition}
\label {a2rdiv}
Suppose that $Z$ is an $\apclass {2}$-regular quasi-normed lattice on $(S \times \Omega, \nu \times \mu)$.
Then lattice $Z^{\frac 1 2} \lclassg {1}^{\frac 1 2}$ is $\apclass {1}$-regular.
\end {proposition}
Indeed, suppose that $f \in Z^{\frac 1 2} \lclassg {1}^{\frac 1 2} = Z^{\frac 1 2} \lclassg {2}$ with norm $1$,
so there exist some $g \in Z$ and $h \in \lclassg {2}$ with norms at most $2$ such that
$g \geqslant 0$ almost everywhere and $f = g^{\frac 1 2} h$.
Let $\weightw$ be a suitable $\apclass {2}$-majorant for $g$ in $Z$.
Then we also have $\weightw^{-1} \in \apclass {2}$ and 
\begin {multline*}
\|M f\|_{Z^{\frac 1 2} \lclassg {2}} \leqslant
\left\|\weightw^{\frac 1 2}\right\|_{Z^{\frac 1 2}} \left\|(M f) \weightw^{- \frac 1 2} \right\|_{\lclassg {2}} =
\left\|\weightw\right\|_{Z}^{\frac 1 2} \left\|(M f) \weightw^{- \frac 1 2} \right\|_{\lclassg {2}} \leqslant
\\
c \left\|(M f) \weightw^{- \frac 1 2} \right\|_{\lclassg {2}} =
c \left\|(M f) \right\|_{\lclass {2} {\weightw^{\frac 1 2}}} \leqslant
\\
c' \left\|f \right\|_{\lclass {2} {\weightw^{\frac 1 2}}} =
c' \left\|f \weightw^{-\frac 1 2} \right\|_{\lclassg {2}} =
c' \left\|\left(\frac g \weightw\right)^{\frac 1 2} h \right\|_{\lclassg {2}} \leqslant 
c' \left\|h \right\|_{\lclassg {2}} \leqslant c''
\end {multline*}
with some suitable constants $c$, $c'$ and $c''$, so $M$ is bounded in the lattice
$Z^{\frac 1 2} \lclassg {1}^{\frac 1 2}$ which is thus $\apclass {1}$-regular as claimed.

Now we can prove Theorem~\ref {themcr2}.
Since $X$ is $2$-convex, we may apply Theorem~\ref {btsb}
to lattice $Y = X^2$ and obtain $\apclass {2}$-regularity of lattice $Y'$ by the assumed $\apclass {2}$-nondegeneracy
of operator $T$.  By Proposition~\ref {aptoconj} lattice $Y^{\frac 1 2} = X$ is then $\apclass {1}$-regular,
and by Proposition~\ref {a2rdiv} lattice
$
Y'^{\frac 1 2} \lclassg {1}^{\frac 1 2} = Y'^{\frac 1 2} \lclassg {\infty}'^{\frac 1 2} =
\left(Y^{\frac 1 2} \lclassg {\infty}^{\frac 1 2}\right)' = X'
$
is also $\apclass {1}$-regular as claimed.

\section {Main lemma revisited}

\label {mlr}

Recall that a lattice $X$ is called $\apclass {p}$-regular if functions from $X$ admit $\apclass {p}$ majorants
with the appropriate control on the norm; see also Definition~\ref {frdef} in Section~\ref {doapar} below.
Lattice $X$ is $\apclass {\infty}$-regular if and only if it is $\apclass {p}$-regular with some $p < \infty$.
$\apclass {1}$-regularity is equivalent to the boundedness of the Hardy-Littlewood
maximal operator (see, e.~g., \cite [Proposition~1] {rutsky2011en}).

The following result was established in~\cite [Theorem~8] {rutsky2014sm}
with the help of a fixed point theorem under an additional assumption that
$X$ is a Banach lattice satisfying the Fatou property.
However, we will now see that
for the proof
it suffices to carry out
a slightly modified version of estimate~\cite [(6)] {rutsky2014sm} with the appropriate majorants.
\begin {proposition}
\label {a1apt}
Suppose that $X$ is a quasi-Banach lattice of measurable functions on $(S \times \Omega, \mu \times \nu)$
such that $X$ is $\apclass {p}$-regular with some $1 \leqslant p < \infty$ and $X^\delta$ is $\apclass {1}$-regular with
some $\delta > 0$.  Then $X$ is $\apclass {1}$-regular with 
an appropriate estimate for the
constants depending only on the corresponding $\apclass {p}$-regularity constants of $X$,
$\apclass {1}$-regularity constants of $X^\delta$ and the value of $\delta$.
\end {proposition}

Indeed, let $f \in X$.  Then there exists an $\apclass {p}$-majorant $\weightw$ for $f$ in $X$,
and in turn there exists an $\apclass {1}$-majorant $\weightu$ for $\weightw^\delta$ in $X^\delta$.
We fix some $\omega \in \Omega$ such that
$\weightw (\cdot, \omega) \in \apclass {p}$ and $\weightu (\cdot, \omega) \in \apclass {1}$,
and let $B (x, r) \subset S$, $x \in S$, $r > 0$, be an arbitrary ball in $S$.
Sequential application of the $\apclass {p}$ condition satisfied by weight $\weightw$,
the Jensen inequality with convex function $t \mapsto t^{-\delta (p - 1)}$, $t > 0$,
and the $\apclass {1}$ condition satisfied by the weight $\weightu$
yield
\begin {multline}
\label {a1condb}
\frac 1 {\nu (B (x, r))} \int_{B (x, r)} |f (\cdot, \omega)| \leqslant
\frac 1 {\nu (B (x, r))} \int_{B (x, r)} \weightw (\cdot, \omega) \leqslant
\\
c \left[ \frac 1 {\nu (B (x, r))} \int_{B (x, r)} [\weightw (\cdot, \omega)]^{-\frac 1 {p - 1}} \right]^{-(p - 1)} =
\\
c \left[ \frac 1 {\nu (B (x, r))} \int_{B (x, r)} [\weightw (\cdot, \omega)]^{-\frac 1 {p - 1}} \right]^{-\delta (p - 1) \cdot \frac 1 \delta} \leqslant
\\
c \left[ \frac 1 {\nu (B (x, r))} \int_{B (x, r)} [\weightw (\cdot, \omega)]^\delta  \right]^{\frac 1 \delta} \leqslant
\\
c \left[ \frac 1 {\nu (B (x, r))} \int_{B (x, r)} \weightu (\cdot, \omega) \right]^{\frac 1 \delta} \leqslant
c' \left[\weightu (x, \omega)\right]^{\frac 1 \delta}
\end {multline}
for almost all $x \in S$ with suitable constants $c$ and $c'$.
Since $\omega$, $x$ and $B$ are arbitrary, \eqref {a1condb} implies that $M f \leqslant c' \weightu^{\frac 1 \delta}$
almost everywhere, so
$\|M f\|_X \leqslant c' \|\weightu\|_{X^\delta}^{\frac 1 \delta} \leqslant
c'' \left\|\weightw^{\delta}\right\|_{X^\delta}^{\frac 1 \delta} =
c'' \left\|\weightw\right\|_{X} \leqslant c''' \|f\|_X$
with some appropriate constants $c''$ and $c'''$.  Thus $M$ is bounded in $X$ with an appropriate estimate
of the norm, and so lattice $X$ is suitably $\apclass {1}$-regular.

\begin {proposition}
\label {ainfainf}
Let $X$ be a normed lattice on $(S \times \Omega, \nu \times \mu)$ such that $X'$ is norming for $X$.
Suppose that both $X$ and $X'$ are $\apclass {\infty}$-regular.
Then both $X$ and $X'$ are $\apclass {1}$-regular.
\end {proposition}
Indeed, since $X$ and $X'$ are $\apclass {\infty}$-regular, they are also $\apclass {p}$-regular with some $p < \infty$.
By Proposition~\ref {aptoconj} both $X'^{\frac 1 p}$ and $X^{\frac 1 p}$ are then $\apclass {1}$-regular,
and it remains to apply Theorem~\ref {a1apt} to $X$ and to $X'$
with $\delta = \frac 1 p$.

\section {Divisibility of $\apclass {p}$-regularity}

\label {doapar}

It is often convenient to think about Muckenhoupt weights in terms of the Jones factorization theorem
(see, e.~g., \cite [Chapter~5, \S 5.3] {stein1993}:
$\weightw \in \apclass {p}$ if and only if $\weightw = \omega_0 \omega_1^{1 - p}$ with some weights
$\omega_0, \omega_1 \in \apclass {1}$ with the appropriate estimates on the constants.
This makes it intuitive that, for example, division by the $\apclass {1}$ weights turns $\apclass {p}$ weights
into $\apclass {p + 1}$ weights, which is the main insight behind the divisibility theorem for $\apclass {p}$-regularity
\cite [Theorem~2] {rutsky2011en}: under certain assumptions on Banach lattices $X$ and $Y$,
if lattice $X Y$ is $\apclass {p}$-regular
and lattice $Y$ is $\apclass {1}$-regular then lattice $X$ is $\apclass {p + 1}$-regular.

However, in the present work a somewhat more general problem arises: we need to make sure that a lattice
$X$ admits majorants $\weightw$ such that $\weightw^{-1} \in \apclass {\infty}$ based on the assumption
that lattice $(X Y)^{\delta}$ is $\apclass {1}$-regular with an $\apclass {1}$-regular lattice $Y$ and some $\delta > 0$.
With that in mind we introduce the following notions; see also \cite [\S 1] {rutsky2015cen}.
\begin {definition}
\label {fdef}
Let $\alpha, \beta \geqslant 0$.
We say that a weight $\weightw$ on $(S \times \Omega, \nu \times \mu)$ belongs to
class $\frclass {\alpha} {\beta}$ with a constant $C$ if there exist two weights $\omega_0, \omega_1 \in \apclass {1}$
with constant $C$ such that $\weightw = \frac {\omega_0^\alpha} {\omega_1^\beta}$.
\end {definition}
\begin {definition}
\label {frdef}
Let $\alpha, \beta \geqslant 0$, and suppose that $X$ is a quasi-normed lattice
on $(S \times \Omega, \nu \times \mu)$.  We say that $X$ is $\frclass {\alpha} {\beta}$-regular
with constants $(C, m)$ if for any $f \in X$ there exists a majorant $\weightw \in X$, $\weightw \geqslant |f|$ such that
$\|\weightw\|_X \leqslant m \|f\|_X$
and $\weightw \in \frclass {\alpha} {\beta}$ with constant $C$.
\end {definition}
``$\frclass {} {}$'' in the notation $\frclass {\alpha} {\beta}$ stands for ``factorizable weight'',
and the properties of the $\apclass {1}$ weights imply that
at least in the local terms $\omega_0$ roughly represents the ``poles'' of the weight $\weightw$
where the weight takes relatively
large values, whereas $\omega_1$ represents the ``zeroes'' of $\weightw$ where the weight is
relatively small.  The corresponding factorization is generally not unique.

Since $\omega \in \apclass {1}$ implies $\omega^\delta \in \apclass {1}$ for $0 < \delta \leqslant 1$,
we see that $\frclass {\alpha} {\beta} \subset \frclass {\alpha_1} {\beta_1}$ for $\alpha \leqslant \alpha_1$
and $\beta \leqslant \beta_1$.  Likewise, $\frclass {\alpha} {\beta}$-regularity of a lattice $X$ implies its
$\frclass {\alpha_1} {\beta_1}$-regularity.

It is easy to see that these properties are closely related to $\apclass {p}$-regularity.
\begin {proposition}
\label {frapr}
Suppose that $\alpha > 0$, $\beta \geqslant 0$ and $\weightw$ is a weight on $(S \times \Omega, \nu \times \mu)$.
Then $\weightw \in \frclass {\alpha} {\beta}$ if and only if $\weightw^{\frac 1 \alpha} \in \apclass {\frac \beta \alpha + 1}$,
and $\weightw \in \frclass {0} {\beta}$ if and only if $\weightw^{-\frac 1 \beta} \in \apclass {1}$
with the appropriate estimates on the constants.
\end {proposition}
Indeed, it suffices to observe that 
$\frac {\omega_0^\alpha} {\omega_1^\beta} = \left(\omega_0 \omega_1^{-\frac \beta \alpha}\right)^\alpha$ for $\alpha > 0$.

Proposition~\ref {frapr} yields at once the corresponding result for $\frclass {\alpha} {\beta}$-regularity.
\begin {proposition}
\label {frapreg}
Let $X$ be a quasi-normed lattice on $(S \times \Omega, \nu \times \mu)$, and suppose that $\alpha > 0$, $\beta \geqslant 0$.
Lattice $X$ is $\frclass {\alpha} {\beta}$-regular if and only if lattice $X^{\frac 1 \alpha}$ is
$\apclass {\frac \beta \alpha + 1}$-regular.
\end {proposition}
Incidentally, as a corollary we get yet another characterization of the property $\log \weightw \in \BMO$ and
the corresponding $\BMO$-regularity in terms of $\weightw \in \frclass {\alpha} {\beta}$ with some $\alpha$ and $\beta$
and, respectively, $\frclass {\alpha} {\beta}$-regularity of lattice $X$.

Notation $\frclass {\alpha} {\beta}$ allows convenient computations for exponents and products of weights.
The following property is immediate from the definitions.
\begin {proposition}
\label {frappower}
Suppose that $\alpha, \beta, \gamma \geqslant 0$ and $\weightw \in \frclass {\alpha} {\beta}$.
Then $\weightw^\gamma \in \frclass {\gamma \alpha} {\gamma \beta}$ with the same constants.
If $\weightw > 0$ almost everywhere then $\weightw^{-\gamma} \in \frclass {\gamma \beta} {\gamma \alpha}$
with the same constants.
If a lattice $X$ is $\frclass {\alpha} {\beta}$-regular and $\gamma > 0$ then lattice $X^\gamma$ is 
$\frclass {\gamma \alpha} {\gamma \beta}$-regular with the same constants.
\end {proposition}
\begin {proposition}
\label {frapmult}
Suppose that $\alpha_0, \alpha_1, \beta_0, \beta_1 \geqslant 0$,
$\weightw_0 \in \frclass {\alpha_0} {\beta_0}$ and $\weightw_1 \in \frclass {\alpha_1} {\beta_1}$.
then $\weightw_0 \weightw_1 \in \frclass {\alpha_0 + \alpha_1} {\beta_0 + \beta_1}$ with the appropriate estimates
on the constants.
Likewise, if $X$ and $Y$ are some lattices on $(S \times \Omega, \nu \times \mu)$
such that $X$ is $\frclass {\alpha_0} {\beta_0}$-regular and $Y$ is $\frclass {\alpha_1} {\beta_1}$-regular
then lattice $X Y$ is $\frclass {\alpha_0 + \alpha_1} {\beta_0 + \beta_1}$-regular
with the appropriate estimates on the constants.
\end {proposition}
Indeed, since the sets of $\apclass {1}$ weights with constant at most $C$
are logarithmically convex (see \eqref {abp1} below), it is easy to see that
if $\weightw_0 = \frac {\omega_{00}^{\alpha_0}} {\omega_{01}^{\beta_0}} \in \frclass {\alpha_0} {\beta_0}$ and
$\weightw_1 = \frac {\omega_{10}^{\alpha_1}} {\omega_{11}^{\beta_1}} \in
\frclass {\alpha_0} {\beta_0}\in \frclass {\alpha_1} {\beta_1}$
with some appropriate $\omega_{jk} \in \apclass {1}$
then
$$
\weightw_0 \weightw_1 =
\frac {\left(\omega_{00}^{\frac {\alpha_0} {\alpha_0 + \alpha_1}} \omega_{10}^{\frac {\alpha_1} {\alpha_0 + \alpha_1}} \right)^{\alpha_0 + \alpha_1}}
 {\left(\omega_{01}^{\frac {\beta_0} {\beta_0 + \beta_1}} \omega_{11}^{\frac {\beta_1} {\beta_0 + \beta_1}} \right)^{\beta_0 + \beta_1}}
\in \frclass {\alpha_0 + \alpha_1}
{\beta_0 + \beta_1}
$$
with an appropriate estimate for the constant.

It is remarkable that the statement of Proposition~\ref {frapmult} can be reversed not only for weights
but also for lattices.
The following result is a generalization of~\cite [Theorem~2] {rutsky2011en};
in the proof of Theorem~\ref {themcr} in Section~\ref {noaro} below it is applied
with $\alpha_1 = 2$, $\alpha_0 = 1$, $\beta_0 = \beta_1 = 0$,

\begin {theorem}
\label {frdiv}
Suppose that $X$ and $Y$ are quasi-Banach lattices on $(S \times \Omega, \nu \times \mu)$ satisfying the
Fatou property, $X Y$ is $\frclass {\alpha_1} {\beta_1}$-regular and $Y$ is $\frclass {\alpha_0} {\beta_0}$-regular.
Then lattice $X$ is $\frclass {\alpha_1 + \beta_0} {\beta_1 + \alpha_0}$-regular.
\end {theorem}
Examining the case of
weighted $\lclass {\infty} {\weightw}$ lattices with suitable weights shows that the conclusion
of Theorem~\ref {frdiv} is sharp in the sense that the indexes of regularity cannot be replaced by smaller values.

A complete proof of theorem~\ref {frdiv} is given in Section~\ref {frdivproof} below.
A weaker statement can be obtained
directly from \cite [Theorem~2] {rutsky2011en}; however, the resulting indexes of regularity are too crude for
our purposes.  However, we may deduce the case needed in the present work
from the following recently obtained result,
which seems to be somewhat less involved technically
than the proof of Theorem~\ref {frdiv} in full generality that, among other
things, uses a fixed point theorem.

\begin {theorem} [{\cite [Theorem~14] {rutsky2015cen}}]
\label {frabduality}
Suppose that $X$ is a Banach lattice of measurable functions on $S \times \Omega$ satisfying the Fatou property
and $\alpha > 1$, $\beta > 0$.  Then $X$ is $\frclass {\alpha} {\beta}$-regular if and only if
$X'$ is $\frclass {\beta + 1} {\alpha - 1}$-regular.
\end {theorem}

Indeed, suppose that under the conditions of Theorem~\ref {frdiv} both lattices $X$ and $Y$ are
$r$-convex with some $r > 0$ such that $\alpha_0 r > 1$ and $(\alpha_1 + \beta_0) r > 1$;
these conditions are satisfied
in the application to the proof of Theorem~\ref {themcr} in Section~\ref {noaro} below
with some sufficiently close to $1$ value of $r > 1$.
Then lattice $(X^r)^{\frac 1 2} (Y^r)^{\frac 1 2}$ is $\frclass {\frac {\alpha_1 r} 2} {\frac {\beta_1 r} 2}$-regular
and lattice $Y^r$ is $\frclass {\alpha_0 r} {\beta_0 r}$-regular
by Proposition~\ref {frappower}, so lattice $(Y^r)'$ is $\frclass {\beta_0 r + 1} {\alpha_0 r - 1}$-regular
by Theorem~\ref {frabduality}, thus lattice $(Y^r)'^{\frac 1 2}$ is
$\frclass {\frac {\beta_0 r + 1} 2} {\frac {\alpha_0 r - 1} 2}$-regular by Proposition~\ref {frappower}.
By the Lozanovsky factorization theorem \cite {lozanovsky1969} we have
$\lclassg {1} = (Y^r) (Y^r)'$,
and lattice $(X^r)^{\frac 1 2} \lclassg {1}^{\frac 1 2} = (X^r)^{\frac 1 2} (Y^r)^{\frac 1 2} (Y^r)'^{\frac 1 2}$
is $\frclass {\frac {(\alpha_1 + \beta_0) r + 1} 2} {\frac {(\beta_1 + \alpha_0) r - 1} 2}$-regular by 
Proposition~\ref {frapmult}, which by Theorem~\ref {frabduality} implies that
lattice
$\left[(X^r)^{\frac 1 2} \lclassg {1}^{\frac 1 2}\right]' = (X^r)'^{\frac 1 2}$ is
$\frclass {\frac {(\beta_1 + \alpha_0) r + 1} 2} {\frac {(\alpha_1 + \beta_0) r - 1} 2}$-regular,
and thus lattice $(X^r)'$ is $\frclass {(\beta_1 + \alpha_0) r + 1} {(\alpha_1 + \beta_0) r - 1}$-regular
by Proposition~\ref {frappower}.  Applying Theorem~\ref {frabduality} to lattice $(X^r)'$ yields
$\frclass {(\alpha_1 + \beta_0) r} {(\beta_1 + \alpha_0) r}$-regularity of lattice $X^r$,
which by Proposition~\ref {frappower} implies the required
$\frclass {\alpha_1 + \beta_0} {\beta_1 + \alpha_0}$-regularity of lattice $X$.

\section {An estimate for nondegenerate operators}

\label {apono}

It is well known that if $T$ is a nondegenerate operator in the sense of Definition~\ref {nondegop}
then the boundedness of $T$ in $\lclass {2} {\weightw^{-\frac 1 2}}$ implies that $\weightw \in \apclass {2}$.
However, in quantitative terms
the standard argument establishing this (see, e.~g., \cite [Chapter~5, \S 4.6] {stein1993})
only yields an estimate
$[\weightw]_{\apclass {2}} \leqslant C \|T\|_{\lclass {2} {\weightw^{-\frac 1 2}} \to \lclass {2} {\weightw^{-\frac 1 2}}}^4$,
which is too rough for the proof of Theorem~\ref {themcr} in Section~\ref {noaro} below to work in full generality.
The value~$[\weightw]_{2}$ cannot be estimated
in terms of $C \|T\|_{\lclass {2} {\weightw^{-\frac 1 2}} \to \lclass {2} {\weightw^{-\frac 1 2}}}$; see \cite [\S 8.B] {hytonenperez2013}.

Nevertheless, securing an additional restriction on the doubling constant of
either the weight $\weightw$ or the weight $\weightw^{-1}$ leads to a suitable estimate.
We denote by $\lambda_n$ the Legbesgue measure on $\mathbb R^n$.
\begin {proposition}
\label {nondegwd}
Suppose that $T$ is a nondegenerate operator that is bounded in $\lclass {2} {\weightw^{-\frac 1 2}}$
with a weight $\weightw$ on $(\mathbb R^n \times \Omega, \lambda_n \times \mu)$
such that either $\weightw$ or $\weightw^{-1}$ satisfies the doubling condition with a constant $c_{\weightw}$.
Then
\begin {equation}
\label {nondegwdeq}
[\weightw]_{\apclass {2}} \leqslant c_T C_{\weightw} \|T\|_{\lclass {2} {\weightw^{-\frac 1 2}} \to \lclass {2} {\weightw^{-\frac 1 2}}}^2
\end {equation}
with a constant $c_T$ independent of the weight $\weightw$ and a constant $C_\weightw$ depending only on $c_\weightw$.
\end {proposition}
Indeed, let $m = \|T\|_{\lclass {2} {\weightw^{-\frac 1 2}} \to \lclass {2} {\weightw^{-\frac 1 2}}}$ under the assupmtions of
Proposition~\ref {nondegwd}.
The argument in~\cite [Proposition~19] {rutsky2011en} shows that
\begin {equation}
\label {singnd}
\int |T f (\cdot)|^2 \weightw (\cdot, \omega) \leqslant 2 m^2 \int |f (\cdot)|^2 \weightw (\cdot, \omega)
\end {equation}
for almost all $\omega \in \Omega$ and all $f \in \lclass {2} {\weightw^{-\frac 1 2} (\cdot, \omega)}$.

Suppose that $B$ is a ball in $\mathbb R^n$ and let $B' = B + r x_0$ with $r > 0$ and $x_0 \in \mathbb R^n$ taken
from the definition of a nondegenerate operator (Definition~\ref {nondegop}) as applied to $T$.
It is easy to see that the boundedness of $T$ implies that both $\weightw (\cdot, \omega)$
and $\weightw^{-1} (\cdot, \omega)$
are locally summable for almost all $\omega \in \Omega$.
Substituting
the condition \eqref {nondegsoeq} from the definition of a nondegenerate operator into \eqref {singnd}, we see that
\begin {multline}
\label {mbe}
2 m^2 \int_B f^2 (\cdot) \weightw (\cdot, \omega) \geqslant
\int |T f (\cdot)|^2 \weightw (\cdot, \omega) \geqslant
\\
\int_{B'} |T f (\cdot)|^2 \weightw (\cdot, \omega) \geqslant
c^2 \left(\frac 1 {|B|} \int_B f\right)^2 \int_{B'} \weightw (\cdot, \omega)
\end {multline}
for almost all $\omega \in \Omega$ and all $f \in \lclass {2} {\weightw^{-\frac 1 2} (\cdot, \omega)}$ such that
$f \geqslant 0$ and
$\supp f \subset B$.
Putting $f = \weightw^{-1} \chi_B$ into \eqref {mbe} yields
\begin {equation}
\label {shiftat}
\left(\frac 1 {|B|}\int_{B'} \weightw (\cdot, \omega)\right) \left(\frac 1 {|B|} \int_B \weightw^{-1} (\cdot, \omega)\right) \leqslant
2 c^{-2} m^2.
\end {equation}
Since the balls $B = B (x, r)$ and $B' = B (x + r x_0, r)$ are mutually comparable in the sense that
$B' \subset B (x, r (1 + |x_0|))$ and $B \subset B (x + r x_0, r (1 + |x_0|))$,
the doubling condition of either the weight $\weightw$ or the weight $\weightw^{-1}$ implies that
one of the factors on the left-hand side of \eqref {shiftat} is suitably comparable to a similar factor
with either $B$ replaced by $B'$ or vice versa. This observation yields~\eqref {nondegwdeq},
since both $B$ and $B'$ are arbitrary balls of $\mathbb R^n$.

We apply Proposition~\ref {nondegwd} to the situation arising in Theorem~\ref {btsb}.
\begin {proposition}
\label {btsbge}
Suppose that $Y$ is a Banach lattice on $(\mathbb R^n \times \Omega, \lambda_n \times \mu)$ with an order continuous norm,
and let $T$ be a nondegenerate operator acting boundedly in $Y^{\frac 1 2}$.
Suppose also that lattice $Y'$ is $\frclass {\alpha} {1}$-regular with some $\alpha > 0$.
Then for every $f \in Y'$
there exists a majorant
$\weightw \geqslant |f|$, $\|\weightw\|_{Y'} \leqslant m_2 \|f\|_{Y'}$, such that
\begin {equation}
\label {btsbgeest}
[\weightw]_{\apclass {2}} \leqslant C_2 \|T\|_{Y^{\frac 1 2} \to Y^{\frac 1 2}}^2
\end {equation}
with some constants $(C_2, m_2)$ independent of $\weightw$ and $\|T\|_{Y^{\frac 1 2} \to Y^{\frac 1 2}}$.
\end {proposition}
To prove Proposition~\ref {btsbge} we need to show that it is possible to take weights $\weightw$ in the conclusion
of Theorem~\ref {btsb} that also satisfy $\weightw \in \frclass {\alpha} {1}$ with a suitable control on the norm.
To do this we adapt the fixed point argument from the proof of \cite [Theorem~8] {rutsky2014sm}.
This requires a few preparations.

We introduce the following sets of Muckenhoupt weights for $p > 1$:
$$
\abp {p} {C} = \left\{ \weightw \in \apclass {p} \mid [\weightw]_{\apclass {p}} \leqslant C \right\},
$$
\begin {equation}
\label {abp1}
\abpo {1} {C} = \left\{ \weightw \in \apclass {1} \mid \esssup \frac {M \weightw} {\weightw} \leqslant C \right\}.
\end {equation}
Here ``$B \apclass {p}$'' denotes ``the ball of $\apclass {p}$'',
and ``$(\mathrm {MC})$'' indicates that these sets are defined by the Muckenhoupt condition to avoid confusion
with earlier work (e.~g. \cite [Section~3] {rutsky2011en}),
where different (for $p > 1$) sets $\abpo {p} {C}$ were used.
The latter have the advantage of being convex and they can also be used to establish the results of the present work;
however, we do not need the convexity, and the basic facts about sets $\abp {p} {C}$ seem to be simpler.
Such a definition is more in line with the rest of the arguments. 

\begin {proposition}
\label {abplccm}
Sets $\abp {p} {C}$ are logarithmically convex and closed with respect to the convergence in measure.
\end {proposition}
Indeed, the logarithmic convexity follows at once from the H\"older inequality, and the closedness
with respect to the convergence in measure is obtained by twice applying the Fatou lemma: if
$\weightw_n \in \abp {p} {C}$ and $\weightw_n \to \weightw$ almost everywhere then
\begin {multline*}
\frac 1 {\nu (B)} \int_B \weightw (\cdot, \omega) \leqslant \liminf_n \frac 1 {\nu (B)} \int_B \weightw_n (\cdot, \omega)
\leqslant
\\
C \liminf_n \left(\frac 1 {\nu (B)} \int_B \weightw_n^{-\frac 1 {p - 1}} (\cdot, \omega) \right)^{-(p - 1)} =
\\
C \left(\limsup_n \frac 1 {\nu (B)} \int_B \weightw_n^{-\frac 1 {p - 1}} (\cdot, \omega) \right)^{-(p - 1)} \leqslant
\\
C \left(\liminf_n \frac 1 {\nu (B)} \int_B \weightw_n^{-\frac 1 {p - 1}} (\cdot, \omega) \right)^{-(p - 1)} \leqslant
\\
C \left(\frac 1 {\nu (B)} \int_B \weightw^{-\frac 1 {p - 1}} (\cdot, \omega) \right)^{-(p - 1)}
\end {multline*}
for all balls $B \subset S$ and almost all $\omega \in \Omega$, so $\weightw \in \abp {p} {C}$.

According to Proposition~\ref {frapr}, we can define for $\alpha > 0$, $\beta \geqslant 0$
the corresponding sets of $\frclass {\alpha} {\beta}$ weights with a control on the constant by
$$
\fbab {\alpha} {\beta} {C} = \left\{\weightw^\alpha \mid \weightw \in \abp {\frac \beta \alpha + 1} {C} \right\},
$$
$$
\fbab {0} {\beta} {C} = \left\{\weightw^{-\beta} \mid \weightw \in \abpo {1} {C}, \text { $\weightw > 0$ almost everywhere} \right\}.
$$
Consequently, these sets are also logarithmically convex and closed with respect to the convergence in measure.

\begin {proposition}
\label {wambspace}
Suppose that $Z$ is a Banach\footnote {It is easy to see that Proposition~\ref {wambspace} also holds true for
quasi-normed lattices $Z$.}
lattice on a $\sigma$-finite measurable space,
$\omega_1 \in Z$, $\omega_1 > 0$ almost everywhere,
$E \subset Z$ is a bounded set in $Z$ such that $h \geqslant \omega_1$ for all $h \in E$.
Then there exists some weight $\omega$, $\omega > 0$ almost everywhere, such that
$D = \{\log \weightw \mid \weightw \in E\}$ is a bounded set in $\lclass {2} {\omega^{-\frac 1 2}}$.
\end {proposition}
To prove Proposition~\ref {wambspace}, take any $a \in Z'$ such that $\|a\|_{Z'} = 1$ and $a > 0$ almost everywhere, any
$\sigma \in \lclassg {1}$ such that $\|\sigma\|_{\lclassg {1}} = 1$ and $\sigma > 0$ almost everywhere, and
define a weight $\omega = a \wedge \sigma (1 - [\log \omega_1]^-)^{-2}$.
Then $\log \weightw \in D$ implies
\begin {multline*}
\int_{\weightw \geqslant 1} |\log \weightw|^2 \omega \leqslant \int_{\weightw \geqslant 1} |\log \weightw|^2 a
=
\\
\int_{\weightw \geqslant 1} 4 \left|\log \left(\weightw^{\frac 1 2}\right)\right|^2 a \leqslant
4 \int \weightw a \leqslant 4 \|\weightw\|_Z \|a\|_{Z'}  \leqslant 4 \|\weightw\|_Z
\end {multline*}
and
\begin {multline*}
\int_{\weightw < 1} |\log \weightw|^2 \omega = \int_{\weightw < 1} (-\log \weightw)^2 \omega \leqslant
\int (-\log \omega_1)^2 \omega =
\\
\int ([\log \omega_1]^-)^2 \omega 
\leqslant
\int ([\log \weightw]^-)^2 \sigma (1 - [\log \omega_1]^-)^{-2}
\leqslant \int \sigma = 1,
\end {multline*}
so indeed $D$ is a bounded set in $\lclass {2} {\omega^{-\frac 1 2}}$.

We now begin the proof of Proposition~\ref {btsbge}.
For convenience, let $X = Y'$; lattice $X$ always has the Fatou property.
Let $C$ be the constant from Theorem~\ref {btsb}.  We introduce a set
\begin {multline*}
B_T = \left\{ \weightw \in X \mid \weightw \geqslant 0, \right.
\\
\left.
\int |T g|^2 \weightw \leqslant
\left(C \|T\|_{Y^{\frac 1 2} \to Y^{\frac 1 2}}\right)^2
\int |g|^2 \weightw \text { for all $g \in \lclass {2} {\weightw^{-\frac 1 2}}$} \right\}.
\end {multline*}
Theorem~\ref {btsb} shows that this set is nonempty.
By the complex interpolation $B_T$ is logarithmically convex.  The closedness of the set $B_T$
with respect to the convergence in measure is verified routinely
(see, e.~g., the proof of~\cite [Proposition~16] {rutsky2011en}):
if $\weightw_n \in B_T$ and $\weightw_n \to \weightw$ almost everywhere then
we put $W = \sup_n \weightw_n$ and see that by the Fatou lemma and the Lebesgue dominated convergence theorem
\begin {multline}
\label {pest2}
\int |T g|^2 \weightw \leqslant
\liminf_{j \to \infty} \int |T g|^2 \weightw_j \leqslant
\left(C \|T\|_{Y^{\frac 1 2} \to Y^{\frac 1 2}}\right)^2 \liminf_{j \to \infty} \int |g|^2 \weightw_j \leqslant
\\
\left(C \|T\|_{Y^{\frac 1 2} \to Y^{\frac 1 2}}\right)^2
\lim_{j \to \infty}
\int \left[ \frac {\weightw_j} {W} \right] |g|^2 W =
\left(C \|T\|_{Y^{\frac 1 2} \to Y^{\frac 1 2}}\right)^2 \int |g|^2 \weightw
\end {multline}
for all $g \in \lclass {2} {W^{-\frac 1 2}}$, so extending
\eqref {pest2} to all $g \in \lclass {2} {\weightw^{-\frac 1 2}}$ by density yields $\weightw \in B_T$.

{}

Suppose that $f \in X$.  We may assume that $\|f\|_X = 1$ and $f > 0$ almost everywhere.
By the assumptions lattice $X$ is $\frclass {\alpha} {1}$-regular with some constants $(C_1, m_1)$.
Let $0 < \beta \leqslant 1$ be a sufficiently small number to be determined later.
We introduce a set $D = \left\{\log \weightw \mid \weightw \in X, \weightw \geqslant \beta f,
\|\weightw\|_X \leqslant 1 \right\}$
and a set-valued map $\Phi : D \times D \to 2^{D \times D}$
by
\begin {multline*}
\Phi ((\log \weightu, \log \weightv)) = \left\{(\log \weightu_1, \log \weightv_1) \in D \times D \mid \right.
\\
\weightu_1 \in X, \weightv_1 \in X, \|\weightu_1\|_X \leqslant 1, \|\weightv_1\|_X \leqslant 1,
\\
\left.
\weightu_1 \in B_T,
\weightv_1 \in \fbab {\alpha} {1} {C_1},
f \vee \weightu \vee \weightv \leqslant A (\weightu_1 \wedge \weightv_1)
\right\}
\end {multline*}
for all $(\log \weightu, \log \weightv) \in D \times D$
with a sufficiently large constant $A$ to be determined in a moment.

Let $(\log \weightu, \log \weightv) \in D \times D$.  Then $\weightw = f \vee \weightu \vee \weightv \in X$
with $\|\weightw\|_X \leqslant 3$.
Applying Theorem~\ref {btsb} to function $\weightw$ yields a majorant $\weightu_2 \in X$,
$\weightu_2 \geqslant \weightw$, $\|\weightu_2\|_X \leqslant 2 \|\weightw\|_X \leqslant 6$
such that $\weightu_2 \in B_T$.  On the other hand, by the $\frclass {\alpha} {1}$-regularity of $X$
there exists some majorant $\weightv_2 \in X$, $\weightv_2 \geqslant \weightw$,
$\|\weightv_2\|_X \leqslant m_1 \|\weightw\|_X \leqslant 3 m_1$ such that $\weightv_2 \in \fbab {\alpha} {1} {C_1}$.
Setting  $\weightu_1 = \frac 1 6 \weightu_2$, $\weightv_1 = \frac 1 {3 m_1} \weightv_2$ and choosing
$A = 6 \vee 3 m_1$ and $\beta = \frac 1 A$ shows that
$(\log \weightu_1, \log \weightv_1) \in \Phi ((\log \weightu_, \log \weightv))$, so $\Phi$ takes nonempty values.

Now it suffices to establish that map $\Phi$ has a fixed point $(\log \weightu, \log \weightv) \in D \times D$,
$\Phi ((\log \weightu, \log \weightv)) \ni (\log \weightu, \log \weightv)$.  If this is the case
then $f \vee \weightu \vee \weightv \leqslant A (\weightu \wedge \weightv)$,
so $\weightw = A (\weightu \vee \weightv)$ is a majorant of $f$ such that
$\|\weightw\|_X \leqslant 2 A$ and
$\weightw$ is pointwise equivalent to both $\weightu$ and $\weightv$ with constant $A$,
which implies that
$\|T\|_{\lclass {2} {\weightw^{-\frac 1 2}} \to \lclass {2} {\weightw^{-\frac 1 2}}} \leqslant A^2 C \|T\|_{Y^{\frac 1 2} \to Y^{\frac 1 2}}$ and
$\weightw \in \frclass {\alpha} {1}$ with a constant depending only on $A$, $C_1$ and $m_1$.
Thus $\weightw^{-1} \in \frclass {1} {\alpha} = \apclass {\alpha + 1}$ by Proposition~\ref {frappower},
and hence $\weightw^{-1}$ is a doubling weight with an estimate for the doubling constant
depending only on $A$, $C_1$ and $m_1$.
Finally, Proposition~\ref {nondegwd} yields the required estimate~\ref {btsbgeest} with a suitable constant $C_2$.

Thus it suffices to verify that $\Phi$ satisfies the assumptions of the Fan--Kakutani fixed point
theorem~\cite {fanky1952}.  We apply Proposition~\ref {wambspace},
which gives a weight $\omega$ such that $D$ is a bounded set in $\lclass {2} {\omega^{-\frac 1 2}}$.
We endow $D$ with the weak topology of $\lclass {2} {\omega^{-\frac 1 2}}$.
$D$ is a convex set that is closed with respect to the convergence in measure,
and hence $D$ is compact.
Likewise, the graph $\Gamma$ of $\Phi$ is a convex set, so it suffices
to show that $\Gamma$ is closed in the strong topology of $D \times D \times D \times D$, which easily follows
from the closedness of $\Gamma$ with respect to the convergence in measure
by the Fatou property of the lattice $X$.
This concludes the proof of Proposition~\ref {btsbge}.

\section {Proof of the main result}

\label {noaro}

We begin by stating a recently developed (see \cite {hytonenperezrela2012})
quantitative estimate for the reverse H\"older inequality
as it applies to $\apclass {p}$-regularity.

The Fujii-Wilson constant of a weight $\weightw \in \apclass {\infty}$ on $S \times \Omega$,
which gives an equavalent definition of the class $\apclass {\infty}$,
is
$$
[\weightw]_{\apclass {\infty}} = \esssup_{\omega \in \Omega}
\sup_{B} 
\frac
{\int_B M [\chi_B \weightw] (\cdot, \omega)}
{\int_B \weightw (\cdot, \omega)},
$$
where the supremum is taken over all balls $B \subset S$.
This constant is dominated by the Muckenhoupt constant $[\weightw]_{\apclass {p}}$ for any $p$
(see, e.~g., \cite [Proposition~2.2] {hytonenperez2013}).  By \cite [Theorem~1.1] {hytonenperezrela2012} any weight
$\weightw \in \apclass {\infty}$ satisfies the reverse H\"older inequality
with all exponents $1 \leqslant r \leqslant 1 + \frac 1 {c [\weightw]_{\apclass {\infty}}}$ for some constant $c$
depending only on the properties of the underlying space $(S, \nu)$, i.~e.
$$
\left(\frac 1 {\nu (B)} \int \weightw^r (\cdot, \omega) \right)^{\frac 1 r} \leqslant
C \frac 1 {\nu (B)} \int \weightw (\cdot, \omega)
$$
for almost all $\omega \in \Omega$ and all balls $B \subset S$ with some constant $C$ independent of $B$.

If $\weightw \in \apclass {p}$ then $\weightw^{-\frac 1 {p - 1}} \in \apclass {p'}$ and
$[\weightw]_{\apclass {p}} = \left[\weightw^{-\frac 1 {p - 1}}\right]_{\apclass {p'}}$, so it is seen immediately
that $\weightw \in \apclass {p}$ implies $\weightw^r \in \apclass {p}$ and
$\left[\weightw^r\right]_{\apclass {p}} \leqslant c_2 [\weightw]_{\apclass {p}}$ for
all $1 \leqslant r \leqslant 1 + \frac 1 {c_1 [\weightw]_{\apclass {p}}}$ with some constants $c_1$ and $c_2$ independent
of $\weightw$.  This implies the following observation.

\begin {proposition}
\label {aregrh}
Suppose that a quasi-Banach lattice $X$ on $(S \times \Omega, \nu \times \mu)$ is $\apclass {\infty}$-regular
with (Fujii-Wilson) constants $(C_{\apclass {\infty}}, m)$, and $X$ is $\apclass {p}$-regular with some $1 \leqslant p < \infty$.
Then $X^r$ is also $\apclass {p}$-regular
for all $1 < r \leqslant 1 + \frac 1 {c \, C_{\apclass {\infty}}}$ with some constant $c$ independent of $C_{\apclass {\infty}}$.
\end {proposition}

\begin {proposition}
\label {l1clp}
Suppose that a normed lattice $X$ on $(S \times \Omega, \nu \times \mu)$ is $\apclass {p}$-regular.
Then lattices $X^{\theta} \lclassg {1}^{1 - \theta}$ are also $\apclass {p}$-regular for all $0 < \theta < 1$.
\end {proposition}
Indeed, let $r > 1$.  We have
$Z = X^\theta \lclassg {1}^{1 - \theta} = \left(X^r\right)^{\frac \theta r} \lclassg {t}^{1 - \frac \theta r}$
with $t = \frac {1 - \frac \theta r} {1 - \theta} > 1$ and $X^r$ is $\apclass {p}$-regular for small
enough values of $r$ by Proposition~\ref {aregrh},
which implies
(by, e.~g., Propositions~\ref {frapmult} and~\ref {frappower})
that $Z$ is also $\apclass {p}$-regular.

We are now ready to prove implication $3 \Rightarrow 1$ of Theorem~\ref {themcr}.
Suppose that under the conditions of Theorem~\ref {themcr} operator $T$ is bounded in~$X$;
we need to show that lattices $X$ and $X'$ are $\apclass {1}$-regular.
By Proposition~\ref {ainfainf} it is sufficient to show that both $X$ and $X'$ are $\apclass {2}$-regular.

The Fatou property together with $p$-convexity and $q$-concavity assumptions on $X$
imply that both $X$ and $X'$ have order continuous
norm (since, for example, $X' = \left(X^p\right)'^{\frac 1 p} \lclassg {1}^{\frac 1 {p'}}$ and the product of
a couple if Banach lattices has order continuous norm if one of the lattices has it),
so $\lclassg {2} \cap X$ is dense in $X$ and
it is easy to see that $T$ is bounded in $X$ if and only if $T^*$ is bounded in $X'$.  Thus by the symmetry
it suffices to prove that $X'$ is $\apclass {2}$-regular.

Let $1 < r, s \leqslant 2$ and $Y = X^r \lclassg {s'} = (X^{r s})^{\frac 1 s} \lclassg {1}^{1 - \frac 1 s}$.
Since $X$ is $p$-convex, $Y$ is a Banach lattice for all $r$ and $s$ satisfying $r s \leqslant p$.
For clarity we may assume that $p \leqslant 2$.
Let us fix $r = \frac {1 + p} 2 < p$;
then $Y$ is $p_1$-convex with
$p_1 = \left(\frac r p + \frac 1 {s'}\right)^{-1} = \left(\frac {p + 1} {2 p} + \frac 1 {s'}\right)^{-1}$, so 
further restricting $s \leqslant \frac {4 p} {3 p + 1}$ yields estimates
$\frac 1 {s'} \leqslant \frac {p - 1} {4 p}$ and
$1 < \frac {4 p} {3 p + 1} = \left(\frac {p + 1} {2 p} + \frac {p - 1} {4 p}\right)^{-1} \leqslant p_1$.
Thus lattice $Y$ is also $p_Y$-convex with $p_Y = \frac {4 p} {3 p + 1}$ for all $s \leqslant \frac {4 p} {3 p + 1} = p_Y$.

We have
$$
Y^{\frac 1 2} = X^{\frac r 2} \lclassg {s'}^{\frac 1 2} =
X^{\frac r 2} \lclassg {1}^{\frac 1 {2 s'}} =
X^{\frac r 2} \left(\lclassg {1}^{\frac 1 {2 s' \left(1 - \frac r 2\right)}}\right)^{1 - \frac r 2} = 
X^{\frac r 2} \left(\lclassg {(2 - r) s'} \right)^{1 - \frac r 2},
$$
and by the complex interpolation we see that
\begin {equation}
\label {tmest}
\|T\|_{Y^{\frac 1 2}} \leqslant \|T\|_X^{\frac r 2} \|T\|_{\lclassg {(2 - r) s'}}^{1 - \frac r 2} \leqslant
c (s')^{1 - \frac r 2}
\end {equation}
with some constant $c$ independent of $s$, since $\|T\|_{\lclassg {t}} = O (t)$ as $t \to \infty$ for a
Calder\'on-Zygmund operator $T$.

A similar computation shows that
\begin {equation}
\label {yinterp}
Y^\varepsilon \lclassg {t}^{1 - \varepsilon} = X^{r \varepsilon} \lclassg {u}^{1 - r \varepsilon}
\end {equation}
if $\frac 1 u (1 - r \varepsilon) = \varepsilon \frac 1 {s'} + (1 - \varepsilon) \frac 1 t$
for some $1 \leqslant u, t < \infty$ and $0 < \varepsilon < 1$.
We choose $\varepsilon = \frac 1 2$ and $t = p_Y'$.  Then \eqref {yinterp} holds true
with $u = \frac {2 - r} {\frac 1 {s'} + \frac 1 t} = \frac {3 - p} {2 \left(\frac 1 {s'} + \frac 1 t \right)}$.
By making $p$ smaller if necessary we may further assume that $t = p_Y' \geqslant \frac 8 {3 - p}$,
and for all $s \leqslant p_Y = t'$, we have $\frac 1 {s'} \leqslant \frac 1 {t}$ and
$2 \leqslant t \frac {3 - p} {4} = \frac {3 - p} {2 \left(\frac 1 t + \frac 1 t\right)}
\leqslant u \leqslant \frac {3 - p} {2 \frac 1 t} = t \frac {3 - p} {2} \leqslant t$.
Thus $T$ is bounded in $\lclassg {u}$
uniformly in $1 \leqslant s \leqslant p_Y$ for the chosen values of $\varepsilon$ and $t$.
The complex interpolation yields
$$
\|T\|_{Y^\varepsilon \lclassg {t}^{1 - \varepsilon}} \leqslant \|T\|_X^{r \varepsilon} \|T\|_{\lclassg {u}}^{1 - r \varepsilon} \leqslant c_1
$$
with a constant $c_1$ independent of $1 \leqslant s \leqslant p_Y$.

Lattice $Y^\varepsilon \lclassg {t}^{1 - \varepsilon} = Y^{\frac 1 2} \lclassg {t}^{\frac 1 2}$ is $p_2$-convex with
$p_2 = \left(\varepsilon \frac 1 {p_Y}  + (1 - \varepsilon) \frac 1 t\right)^{-1} =
2 \left(\frac 1 {p_Y} + 1 - \frac 1 {t'}\right)^{-1} = 2$,
so
we may apply Theorem~\ref {themcr2} to it.
This shows that $\left(Y^{\frac 1 2} \lclassg {t}^{\frac 1 2}\right)' = Y'^{\frac 1 2} \lclassg {t'}^{\frac 1 2}$ is
$\apclass {1}$-regular, or $\frclass {1} {0}$-regular in terms of Definition~\ref {frdef}.
By Proposition~\ref {frappower}
lattice $Y' \lclassg {t'} = \left(Y'^{\frac 1 2} \lclassg {t'}^{\frac 1 2}\right)^2$ is $\frclass {2} {0}$-regular.
Lattice $\lclassg {t'}$ is $\apclass {1}$-regular, or $\frclass {1} {0}$-regular.  Therefore
by Theorem~\ref {frdiv} lattice $Y'$ is $\frclass {2} {1}$-regular.

Thus \eqref {tmest} by Proposition~\ref {btsbge} implies that lattice $Y'$ is $\apclass {2}$-regular
with constants $(C_3, m_3)$ satisfying
$C_3 \leqslant c_3 \|T\|_{Y^{\frac 1 2} \to Y^{\frac 1 2}}^2 \leqslant c_3 c^2 (s')^{2 - r}$ for some $c_3$ and $m_3$ independent
of $s$.  By Proposition~\ref {aregrh} lattice $(Y')^\rho$ is then $\apclass {2}$-regular for all
$1 \leqslant \rho \leqslant 1 + \frac 1 {c_4 (s')^{2 - r}}$ with a constant $c_4$ independent of $s$.

Observe that $Y' = \left[(X^{r s})^{\frac 1 s} \lclassg {1}^{1 - \frac 1 s}\right]' = (X^{r s})'^{\frac 1 s}$
and $(Y')^\rho = (X^{r s})'^{\frac \rho s}$.
Setting $\rho = 1 + \frac 1 {c_4 (s')^{2 - r}}$ yields
$$
\frac \rho s = \frac {c_4 (s')^{2 - r} + 1} {c_4 (s')^{2 - r}} \cdot \frac {s' - 1}  {s'} =
\frac {c_4 (s')^{3 - r} - c_4 (s')^{2 - r} + s' - 1} {c_4 (s')^{3 - r}}.
$$
Since $0 < 2 - r < 1$, we have $c_4 (s')^{2 - r} \leqslant s' - 1$ for sufficiently large values of $s'$,
that is, for sufficiently close to $1$ values of $s$, so we have $\frac \rho s > 1$ and $\frac s \rho < 1$
for small enough values of $s$.  We fix such an $s$.
Lattice $\left[(Y')^{\rho}\right]^{\frac s \rho} = (X^{r s})'$ is $\apclass {2}$-regular, and
by Proposition~\ref {l1clp} lattice
$X' = \left[\left(X^{r s}\right)^{\frac 1 {r s}}\right]' = (X^{r s})'^{\frac 1 {r s}} \lclassg {1}^{1 - \frac 1 {r s}}$
is also $\apclass {2}$-regular.  This concludes the proof of Theorem~\ref {themcr}.

\section {Proof of Theorem~\ref {frdiv}}

\label {frdivproof}

Compared to \cite [Theorem~2] {rutsky2011en},
the proof of Theorem~\ref {frdiv} essentially
requires only minor technical adjustments; however, to avoid confusion
we provide a complete version of it.  The only apparent difficulty that arises in direct translation of the proof
is that the sets of the corresponding
$\frclass {\alpha} {\beta}$-majorants seem to lack convexity for $\alpha \neq 1$;
however, they are still logarithmically convex,
which suffices to establish closedness of the graph of the map using the same method.
We also use a different ambient space for the map, which makes approximating the problem by restricting the conditions
to sets of finite measure unnecessary.
This modification also
allows us to avoid using a compactness-type result for sets closed with respect to the convergence in measure,
since the standard weak compactness of sets in a weighted $\lclassg {2}$ space suffices.

See Section~\ref {apono} above for the definition of sets $\fbab {\alpha} {\beta} {C}$.

\begin {proposition}
\label {fbmax}
Suppose that $\weightu, \weightv \in \fbab {\alpha} {\beta} {C}$ with a constant $C$.
Then $\weightu \vee \weightv \in \fbab {\alpha} {\beta} {2 C}$.
If $\alpha = 0$ then $\weightu \vee \weightv \in \fbab {\alpha} {\beta} {C}$.
\end {proposition}
Indeed, according to Proposition~\ref {frapr},
if $\alpha > 0$ it suffices to prove Proposition~\ref {fbmax} for the corresponding sets $\abp {p} {C}$
in place of $\fbab {\alpha} {\beta} {C}$.
We have
\begin {multline*}
\frac 1 {\nu (B)} \int_B (\weightu \vee \weightv) (\cdot, \omega) \leqslant
\frac 1 {\nu (B)} \int_B \weightu (\cdot, \omega) + \frac 1 {\nu (B)} \int_B \weightv (\cdot, \omega)\leqslant
\\
C \left[\left(\frac 1 {\nu (B)} \int_B \weightu^{-\frac 1 {p - 1}} (\cdot, \omega)\right)^{-(p - 1)} +
\left(\frac 1 {\nu (B)} \int_B \weightv^{-\frac 1 {p - 1}} (\cdot, \omega)\right)^{-(p - 1)}
\right] \leqslant
\\
2 C \left(\frac 1 {\nu (B)} \int_B (\weightu \vee \weightv)^{-\frac 1 {p - 1}} (\cdot, \omega)\right)^{-(p - 1)}
\end {multline*}
for all balls $B \subset S$ and almost all $\omega \in \Omega$, so indeed
$\weightu \vee \weightv \in \abp {p} {C}$.

In the case $\alpha = 0$ it suffices to show that $\weightu, \weightv \in \abpo {1} {C}$ implies
$\weightu \wedge \weightv \in \abpo {1} {C}$, which follows at once from the estimates
$M (\weightu \wedge \weightv) \leqslant M \weightu \leqslant C \weightu$ and 
$M (\weightu \wedge \weightv) \leqslant M \weightv \leqslant C \weightv$.

We begin the proof of Theorem~\ref {frdiv}.
First of all, since for all $\delta > 0$ the statement of Theorem~\ref {frdiv} for lattices $X$ and $Y$
is equivalent to the same statement for lattices $X^\delta$ and $Y^\delta$ with all indices multiplied by $\delta$,
and since for any quasi-Banach lattice $Z$ lattice $Z^\delta$ is (up to a renorming)
Banach for small enough values of $\delta$ (see, e.~g., \cite [Theorem~3.2.1] {rolewicz1985}),
we may assume that lattices $X Y$, $X$ and $Y$ are all Banach.

Suppose that lattice $X Y$ is $\frclass {\alpha_1} {\beta_1}$-regular with constants $\left(C_{XY}, m_{XY}\right)$
and $Y$ is $\frclass {\alpha_0} {\beta_0}$-regular with constants $(C_Y, m_Y)$.
We can choose $C$ large enough (depending on $C_{XY}$ and $C_Y$) that the 
$\frclass {\alpha_1} {\beta_1}$-majorants in $X Y$ lie in $\fbab {\alpha_1} {\beta_1} {C}$
and the $\frclass {\alpha_0} {\beta_0}$-majorants in $Y$ belong to $\fbab {\alpha_0} {\beta_0} {C}$.

Take any $\omega_0 \in Y$ such that $\|\omega_0\|_Y > 0$.  There exists an $\frclass {\alpha} {\beta}$ majorant
$\omega_1 \in \fbab {\alpha} {\beta} {C}$ for $\omega_0$.
We may assume that $\|\omega_1\|_Y = 1$.  Let
$$
D = \left\{\log \weightw \mid \weightw \in \fbab {\alpha} {\beta} {2 C},
\weightw \geqslant \omega_1, \|\weightw\|_Y \leqslant 2\right\}.
$$

Suppose that $f \in X$; we need to prove that there exists a suitable
$\frclass {\alpha_1 + \beta_0} {\beta_1 + \alpha_0}$-majorant for $f$.  We may assume that $f > 0$ almost everywhere
and that $\|f\|_X = 1$.

Take any function $\log \weightw \in D$.  Then $f \weightw \in X Y$ with norm at most $2$, and
there exist some majorants $g \geqslant f$, $g \in \fbab {\alpha_1} {\beta_1} {C}$,
$\|g\|_{X Y} \leqslant 2 m_{X Y}$.
It is easy to see that (see, e.~g., \cite [(16)] {rutsky2011en})
$$
\|g\|_{X Y} \geqslant \left(1 + \|\omega_0\|_Y\right)^{-1}
\inf_{\stackrel {\|b\|_Y \leqslant 1 + \|\omega_0\|_Y,}{b \geqslant \omega_0}} \left\|g b^{-1}\right\|_X,
$$
so there exists some $b \in Y$, $b \geqslant \omega_0$, $\|b\|_Y \leqslant 2$ such that
$\left\|g b^{-1}\right\|_X \leqslant 4 m_{X Y}$.  Now let $\weightv \geqslant b$,
$\weightv \in \fbab {\alpha_0} {\beta_0} {C}$, $\|\weightv\|_Y \leqslant 2 m_{Y}$ be
an $\frclass {\alpha_0} {\beta_0}$-majorant for~$b$, and let
$\weightw_1 = \left(\frac 1 {2 m_{Y}} \weightv\right) \vee \omega_1$.
Then $\left\|g \weightw_1^{-1}\right\|_X \leqslant
2 m_{Y} \left\|g \weightv^{-1}\right\|_X \leqslant 2 m_{Y} \left\|g b^{-1}\right\|_X \leqslant 8 m_{Y} m_{X Y}$,
$\|\weightw_1\|_{Y} \leqslant 2$ and
$\weightw_1 \in \fbab {\alpha_0} {\beta_0} {2 C}$ by Proposition~\ref {fbmax}.
This shows that a set-valued map $\Phi : D \to 2^D$
defined by
$$
\Phi (\log \weightw) = \left\{ \log \weightw_1 \in D \mid g \geqslant f \weightw, g \in \fbab {\alpha_1} {\beta_1} {C},
\left\|g \weightw_1^{-1}\right\|_X \leqslant 8 m_{Y} m_{X Y}\right\}
$$
takes nonempty values.

If map $\Phi$ has a fixed point $\log \weightw \in D$, $\Phi (\log \weightw) \ni \log \weightw$ then
there exists some function $g \geqslant f \weightw$, $g \in \fbab {\alpha_1} {\beta_1} {C}$ such that
$\left\|g \weightw^{-1}\right\|_X \leqslant 8 m_{Y} m_{X Y}$
and $f_1 = g \weightw^{-1} \in \frclass {\alpha_1 + \beta_0} {\beta_1 + \alpha_0}$ with a suitable estimate of the constant
by Propositions~\ref {frappower} and~\ref {frapmult}, so $f_1$ is then a suitable
$\frclass {\alpha_1 + \beta_0} {\beta_1 + \alpha_0}$-majorant for $f$.

Thus it suffices to show that $\Phi$ satisfies the conditions of the Fan--Kakutani fixed point theorem
\cite {fanky1952}: that $D$ is a compact set in a locally convex linear topological space such that
$\Phi$ has closed graph and that $\Phi$ takes convex closed values that are compact.

By Proposition~\ref {wambspace} there exists a weight $\omega$ such that
$D$ is a bounded set in $\lclass {2} {\omega^{-\frac 1 2}}$.  We endow $D$ with the weak topology
of $\lclass {2} {\omega^{-\frac 1 2}}$.
Since $D$ is convex and closed with respect to the convergence in measure,
$D$ is a closed and bounded convex set in $\lclass {2} {\omega^{-\frac 1 2}}$; hence
$D$ is a compact set.

It is easy to see that the graph $\Gamma$
of $\Phi$ is a convex set, so it suffices to show that $\Gamma$ is a closed set in the strong topology
of the ambient space $\lclass {2} {\omega^{-\frac 1 2}}$.

Suppose that $\log a_j, \log \weightu_j \in D$, $\log a_j \in \Phi (\log u_j)$,
$\log a_j \to \log A \in D$ and $\log \weightu_j \to \log U \in D$ in $\lclass {2} {\omega^{-\frac 1 2}}$;
we need to verify that $\log A \in \Phi (\log U)$.

By passing to a subsequence we may assume that we also have $\log a_j \to \log A$ and $\log \weightu_j \to \log U$
in the sense of the convergence almost everywhere.
We form
a nonincreasing sequence $\log \alpha_j = \bigvee_{k \geqslant j} \log a_k \geqslant \log a_j$
and a nondecreasing sequence $\log \eta_j = \bigwedge_{k \geqslant j} \log \weightu_k \leqslant \log \weightu_j$
of measurable functions
such that $\log \alpha_j \to \log A$ and $\log \eta_j \to \log U$ almost everywhere.

Condition $\log a_j \in \Phi (\log u_j)$ implies that sets
\begin {multline*}
W_j = \left\{ \log g \mid g \geqslant f \eta_j, g \in \fbab {\alpha_1} {\beta_1} {C},
\left\|g \alpha_j^{-1}\right\|_X \leqslant 8 m_{Y} m_{X Y} \right\} \supset
\\
\left\{ \log g \mid g \geqslant f u_j, g \in \fbab {\alpha_1} {\beta_1} {C},
\left\|g a_j^{-1}\right\|_X \leqslant 8 m_{Y} m_{X Y} \right\}
\end {multline*}
are nonempty,
and $W_j$ is a nonincreasing sequence of sets.
Since for all $\log g \in W_1$ we have $g \geqslant f \omega_1 > 0$ almost everywhere and functions
$g$ are uniformly bounded in the weighted Banach lattice $X (\alpha_1)$,
by Proposition~\ref {wambspace} there exists a weight $\omega_2$ such that 
$W_1$ is a bounded set in $\lclass {2} {\omega_2^{-\frac 1 2}}$.  It is easy to see that the sets $W_j$ are
convex and closed with respect to the convergence in measure
(and thus also in the strong topology of lattices satisfying the Fatou property),
so they are compact
in the weak topology of $\lclass {2} {\omega_2^{-\frac 1 2}}$.  This implies that
the set $\bigcap_j W_j$ is nonempty, and so there exists some function $g \in \fbab {\alpha_1} {\beta_1} {C}$
such that $g \geqslant f \eta_j$ and $\left\|g \alpha_j^{-1}\right\|_X \leqslant 8 m_{Y} m_{X Y}$ for all $j$.
Thus $g \geqslant f \left(\bigvee_j \eta_j\right) = f U$ and
$\left\|g a^{-1}\right\|_X \leqslant \bigvee_j \left\|g \alpha_j^{-1}\right\|_X \leqslant 8 m_{Y} m_{X Y}$
by the Fatou property.  The existence of such a function $g$ implies that $\log A \in \Phi (\log U)$ as claimed,
which concludes the proof of Theorem~\ref {frdiv}.

\normalsize
\baselineskip=17pt

\bibliographystyle {plain}

\bibliography {bmora}

\end{document}